\documentclass[12pt, reqno]{amsart}
\usepackage{amsmath, amsthm, amscd, amsfonts, amssymb, graphicx, color, mathrsfs}
\usepackage[bookmarksnumbered, colorlinks, plainpages]{hyperref}
\usepackage[all]{xy}
\usepackage{slashed}

\usepackage{soul}
\usepackage{cancel}
\usepackage{ulem}

\newtheorem{theorem}{Theorem}[section]

\newtheorem{corollary}[theorem]{Corollary}
\theoremstyle{definition}

\newtheorem{conjecture}[theorem]{Conjecture}

\theoremstyle{remark}
\newtheorem{remark}[theorem]{Remark}
\numberwithin{equation}{section}

\begin{document}
\setcounter{page}{1}

\title[ The Stein restriction  problem  on the torus  ]{A note on the Stein restriction conjecture and the restriction problem  on the torus }

\author[D. Cardona]{Duv\'an Cardona}
\address{
  Duv\'an Cardona:
  \endgraf
  Department of Mathematics
  \endgraf
  Pontificia Universidad Javeriana.
  \endgraf
  Bogot\'a
  \endgraf
  Colombia
  \endgraf
  {\it E-mail address} {\rm duvanc306@gmail.com}
  }



\begin{abstract} In this note revision we discuss the Stein restriction problem on arbitrary $n$-torus, $n\geq 2$. In contrast with the usual cases of the sphere, the parabola and the cone, we provide necessary and sufficient conditions  on the Lebesgue indices, by finding conditions which are independent of the dimension $n$. 

\end{abstract} \maketitle

\section{Introduction}
This note is devoted to  the Stein restriction problem on the torus $\mathbb{T}^n,$ $n\geq 2.$ In harmonic analysis, the Stein restriction problem for a smooth hypersurface  $S\subset \mathbb{R}^n,$ asks for the conditions on $p$ and $q,$ $1\leq p,q<\infty,$ satisfying 
\begin{equation}\label{Stein}
   \Vert \hat{f}|_{S}\Vert_{L^q({S},d\sigma)}:= \left(\int\limits_{S}|\hat{f}(\omega)|^qd\sigma(\omega)\right)^{\frac{1}{q}}\leq C\Vert f\Vert_{L^p(\mathbb{R}^n)},
\end{equation} where $d\sigma$ is a surface measure associated to $S,$ the constant $C>0$ is independent of $f,$ and $\widehat{f}|_{S}$ denotes the Fourier restriction of $f$ to $S,$ where
\begin{equation}
    \hat{f}(\xi)=\int\limits_{\mathbb{R}^n}e^{-i2\pi x\cdot \xi} f(x)dx,
\end{equation} is the Fourier transform of $f.$ Let us note that for $p=1,$ the Riemann-Lebesgue theorem implies that $\hat{f}$ is a continuous function on $\mathbb{R}^n$ and we can restrict $\hat{f}$ to every subset $S\subset \mathbb{R}^n.$ On the other hand, if  $f\in L^2(\mathbb{R}^n),$ the Plancherel theorem gives $\Vert f \Vert_{L^2(\mathbb{R}^n)}=\Vert \hat{f} \Vert_{L^2(\mathbb{R}^n)}$ and the Stein restriction problem is trivial by considering that every hypersurface is a subset in $\mathbb{R}^n$ with vanishing Lebesgue measure. So, for $1<p<2,$ a general problem is to  find those hypersurfaces  $S,$ where the Stein restriction problem has sense. However, the central problem in the restriction theory is the following conjecture (due to Stein). It is of particular interest because it is related to Bochner-Riesz multipliers and the Kakeya conjecture.
\begin{conjecture}\label{ConjectureofStein}
Let $S=\mathbb{S}^{n-1}=\{x\in \mathbb{R}^n:|x|=1\}$ be the $(n-1)$-sphere and let $d\sigma$ be the corresponding surface measure. Then \eqref{Stein} holds true if and only if $1\leq p<\frac{2n}{n+1}$ and $q\leq p'\cdot \frac{n-1}{n+1},$ where $p'=p/p-1.$
\end{conjecture}
That the inequalities $1\leq p<\frac{2n}{n+1}$ and $q\leq p'\cdot \frac{n-1}{n+1},$ are necessary conditions for Conjecture \ref{ConjectureofStein} is a well known fact. In  this setting, a celebrated result by Tomas and Stein (see e.g. Tomas \cite{Tomas}) shows that 
\begin{equation}
    \Vert \hat{f}|_{\mathbb{S}^{n-1}} \Vert_{L^2(\mathbb{S}^{n-1},d\sigma)}\leq C_{p,n}\Vert f\Vert_{L^{p}(\mathbb{R}^n)}
\end{equation}
holds true for every $1\leq p\leq \frac{2n+2}{n+3}.$ Surprisingly, a theorem due to Bourgain shows that the Stein restriction conjecture is true for $1<p<p_n$ where $p_n$ is defined  inductively  and $\frac{2n+2}{n+3}<p_n<\frac{2n}{n+1}.$ For instance,
$p(3) = 31/23$. We refer the reader to Tao \cite{Tao2003} for a good introduction and some advances to the restriction theory.

In this paper we will consider the $n$-dimensional torus $\mathbb{T}^n=(\mathbb{S}^1)^n$ modelled on $\mathbb{R}^{2n},$ this means that
\begin{equation}\label{DefiCliffordTorues}
    \mathbb{T}^n=\{(x_{1,1},x_{1,2},x_{2,1},x_{2,2},\cdots,x_{n,1},x_{n,2}):x_{\ell,1}^2+x_{\ell,2}^2=1,\,1\leq \ell\leq n\}.
\end{equation} In this case  $$\mathbb{T}^n\subset {\sqrt{n}}\,\mathbb{S}^{2n-1} \subset \mathbb{R}^{2n}.$$
In order to illustrate our results, we will discuss the case $n=2,$ where
\begin{equation}\label{DefiCliffordTorues2} \mathbb{T}^2\subset {\sqrt{2}}\,\mathbb{S}^{3} \subset \mathbb{R}^{4}.\end{equation}
As it is well known, the $n$-dimensional torus can be understood of different ways. Topologically, $\mathbb{T}^n\sim \mathbb{S}^1\times \cdots \times \mathbb{S}^1,$ where the circle $\mathbb{S}^1$ can be identified with the unit interval $[0,1),$ where we have identified $0\sim 1.$ The case $n=2,$ implies that $\mathbb{T}^2\sim \mathbb{S}^1\times \mathbb{S}^1.$ From differential geometry, a stereographic  projection $\pi$ from $\mathbb{S}^3\setminus \{N\}$ into $\mathbb{R}^3$ gives the following embedding of $(1/\sqrt{2})\mathbb{T}^2\subset \mathbb{S}^3,$
\begin{equation}\label{anothertorus}
    \dot{\mathbb{T}}^2=\{((\sqrt{2}+\cos(\phi))\cos(\theta),(\sqrt{2}+\cos(\phi))\sin(\theta),\sin(\phi)\in \mathbb{R}^3: 0\leq \theta,\phi<2\pi  \},
\end{equation} of the 2-torus in $\mathbb{R}^3.$ At the same time, the Fourier analysis and the geometry on the  torus can be understood  better   by the description of the torus given in \eqref{DefiCliffordTorues}. So, we will investigate the restriction problem on the torus by using \eqref{DefiCliffordTorues2} instead of \eqref{anothertorus}.
In this case, the Stein restriction conjecture for $S=\mathbb{S}^3$ assures that \eqref{Stein} holds true for every $1\leq p<\frac{8}{5}$ and $q\leq \frac{3}{5}p'.$ However, we will prove the following result, where we characterise the Stein restriction problem on $\mathbb{T}^2$.
\begin{theorem}\label{ThrcardonaRest2018}
Let $f\in L^{p}(\mathbb{R}^4).$ Then there exists $C>0,$ independent of $f$ and satisfying the estimate
\begin{equation}\label{cardonaRest2018}
   \Vert \hat{f}|_{\mathbb{T}^2}\Vert_{L^q({\mathbb{T}^2},d\sigma)}:= \left(\int\limits_{\mathbb{T}^2}|\hat{f}(\xi_1,\xi_2,\eta_1,\eta_2)|^qd\sigma(\xi_1,\xi_2,\eta_1,\eta_2)\right)^{\frac{1}{q}}\leq C\Vert f\Vert_{L^{p}(\mathbb{R}^4)},
\end{equation} if and only if $1\leq p< \frac{4}{3}$ and $q\leq p'/3.$  Here, $d\sigma(\xi_1,\xi_2,\eta_1,\eta_2)$ is the usual surface measure associated to $\mathbb{T}^2.$ 
\end{theorem}
An important difference between the restriction problem on the $n$-torus, $n\geq 2,$ and  the Stein-restriction conjecture come from the curvature notion. For example, the  sphere $\mathbb{S}^2$, has Gaussian curvature non-vanishing, in contrast with the 2-torus $\mathbb{T}^2$ where the Gaussian curvature vanishes identically.
In the general case, let us observe that the Stein conjecture for $S=\mathbb{S}^{2n-1}$ asserts that \eqref{Stein} holds true for all $1\leq p<\frac{ 4n}{2n+1}$ and $q\leq \frac{2n-1}{2n+1}p'.$ Curiously, the situation for the $n$-dimensional torus is very different,  as we will see in the following theorem.
\begin{theorem}\label{ThrcardonaRestTn2018}
Let $f\in L^{p}(\mathbb{R}^{2n}),$ $n\geq 2.$ Then there exists $C>0,$ independent of $f$ and satisfying
\begin{equation}\label{cardonaRestTn2018}
   \Vert \hat{f}|_{\mathbb{T}^{n}}\Vert_{L^q({\mathbb{T}^{n}},d\sigma_n)}\leq C_n\Vert f\Vert_{L^{p}(\mathbb{R}^{2n})},
\end{equation} if and only if $1\leq p< \frac{4}{3}$ and $q\leq p'/3.$  Here, $d\sigma_n$ is the usual surface measure associated to $\mathbb{T}^{n}.$ 
\end{theorem}
\begin{remark}
By a duality argument we conclude the following fact: if
 $F\in L^{q'}(\mathbb{T}^{n},d\sigma_n),$ then there exists $C>0,$ independent of $F$ and satisfying
\begin{equation}
   \Vert (Fd\sigma_n)^{\vee} \Vert_{L^{p'}(\mathbb{R}^{2n})}=\left\Vert \int\limits_{ \mathbb{T}^n } e^{i2\pi x\cdot \xi}F(\xi)  d\sigma_n(\xi)   \right\Vert_{L^{p'}({\mathbb{R}^{2n}})}\leq C_n\Vert F\Vert_{L^{q'}{(\mathbb{T}^n,d\sigma_n)}},
\end{equation} if and only if $p'>4$ and $q'\geq (p'/3)'.$ We have denoted by $(Fd\sigma_n)^{\vee}$ the inverse Fourier transform of the measure $\mu:=Fd\sigma_n.$
\end{remark}

We end this introduction by summarising the progress on the restriction conjecture as follows. Indeed, we refer the reader to,

\begin{itemize}
    \item Fefferman \cite{Fefferman1970} and 
 Zygmund \cite{Zygmund74} for the proof of the restriction conjecture in the case $n=2$ (which is \eqref{cardonaRestTn2018} for $n=1$).
 \item   Stein \cite{Stein1986}, Tomas \cite{Tomas} and  Strichartz \cite{Strichartz1977}, for the restriction problem in higher dimensions, with  sharp $(L^q,L^2)$ results
 for hypersurfaces with nonvanishing
Gaussian curvature.  
Some more general classes of surfaces were treated by A. Greenleaf \cite{Greenlaf1981}.
\item Bourgain \cite{Bourgain1991,Bourgain1995b}, Wolff \cite{Wolff1995},  Moyua, Vargas, Vega and Tao \cite{Moyua1996,Moyua1999,Tao1998} who established the so-called bilinear approach.
\item 
Bourgain and Guth \cite{BoutgainGuth2011},   Bennett, Carbery
and Tao \cite{BeCarTao2006}, by the  progress on the case of
nonvanishing curvature, by making use of multilinear restriction estimates. 
\item Finally, Buschenhenke, M\"uller and  Vargas  \cite{Muller}, for a complete list of references as well as the progress on the restriction theory on surfaces of  finite type.
\end{itemize}
The main goal of this note is to give a simple proof of the restriction problem on the torus. This work is organised as follows. In Section \ref{Sec2} we prove Theorem \ref{ThrcardonaRest2018}. We end this note with the proof of Theorem \ref{ThrcardonaRestTn2018}. Sometimes we will use $(\mathscr{F}f)$ for the 2-dimensional Fourier transform of $f$ and $(\mathscr{F}_{\mathbb{R}^n}u)$ for the Fourier transform of a function $u$ defined on $\mathbb{R}^n.$

\section{Proof of Theorem \ref{ThrcardonaRest2018}}\label{Sec2}
In this note we will use the standard notation used for the Fourier analysis on $\mathbb{R}^n$ and the torus (see e.g. Ruzhansky and Turunen \cite{Ruz}).
Throughout this section we will consider the 2-torus $\mathbb{T}^2,$
\begin{equation}
    \mathbb{T}^2=\{(x_{1},x_{2},y_{1},y_{2}):x_{1}^2+x_{2}^2=1,\,y_{1}^2+y_{2}^2=1 \}=\mathbb{S}^1_{(x_1,x_2)}\times \mathbb{S}^1_{(y_1,y_2)}\subset \mathbb{R}^4.
\end{equation} Here,  $\mathbb{T}^2$ will be endowed
with the surface measure $$d\sigma(\xi_1,\xi_2,\eta_1,\eta_2)=d\sigma(\xi_1,\xi_2)d\sigma(\eta_1,\eta_2),$$
where $d\sigma(\xi_1,\xi_2)$ is the usual  `surface measure' defined on $\mathbb{S}^1.$ Indeed, if $(\xi_1,\xi_2)\equiv (\xi_1(\varkappa),\xi_2(\varkappa))=(\cos(2\pi \varkappa),\sin(2\pi\varkappa)),$ $0\leq \varkappa< 1 ,$ then $d\varkappa=d\sigma(\xi_1,\xi_2).$

Conjecture \ref{ConjectureofStein} has been proved by Fefferman for $n=2,$ the corresponding announcement is the following (see Fefferman \cite{Fefferman1970} and Zygmund \cite{Zygmund1974}). \begin{theorem}[Fefferman restriction Theorem]\label{FeffermanStein}
Let $S=\mathbb{S}^{1}=\{x\in \mathbb{R}^2:|x|=1\}$ be the $1$-sphere and let $d\sigma$ be the corresponding `surface measure'. Then \eqref{Stein} holds true if and only if $1\leq p<\frac{4}{3}$ and $q\leq p'/3,$ where $p'=p/p-1.$
\end{theorem}

In order to prove Theorem \ref{ThrcardonaRest2018}, let us consider $1\leq p<\frac{4}{3}$, $q\leq p'/3$ and  $f\in L^{p}(\mathbb{R}^4) .$ By the argument of density we can assume that $f\in C^\infty_{c}(\mathbb{R}^4).$ If $(\xi_1,\xi_2,\eta_1,\eta_2)\in \mathbb{T}^2,$ then
\begin{equation}
    \widehat{f}(\xi_1,\xi_2,\eta_1,\eta_2)=\int\limits_{\mathbb{R}^4}e^{-i2\pi (x\cdot\xi+y\cdot \eta)}f(x,y)dy\,dx,\,\,x=(x_1,x_2),\,y=(y_1,y_2).
\end{equation}
By the Fubini theorem we can write
\begin{align*}
    \widehat{f}(\xi_1,\xi_2,\eta_1,\eta_2)
    =\int\limits_{\mathbb{R}^2}e^{-i2\pi x\cdot \xi}(\mathscr{F}_{y\rightarrow \eta}{f}(x,\cdot))(\eta)dx,\,\,\eta=(\eta_1,\eta_2),
\end{align*} where $(\mathscr{F}_{y\rightarrow \eta}{f}(x,\cdot))(\eta)=\widehat{f}(x,\eta)$ is the 2-dimensional Fourier transform of the function $f(x,\cdot),$ for every $x\in \mathbb{R}^2.$ By writing
\begin{equation}
    \widehat{f}(\xi_1,\xi_2,\eta_1,\eta_2)=\mathscr{F}_{x\rightarrow \xi}(\mathscr{F}_{y\rightarrow \eta}{f}(x,\cdot))(\eta))(\xi),
\end{equation}
for  $1\leq p<\frac{4}{3}$ and $q\leq p'/3,$ the Fefferman restriction theorem gives,
\begin{equation}
    \Vert  \widehat{f}(\xi_1,\xi_2,\eta_1,\eta_2)\Vert_{L^q(\mathbb{S}^1,d\sigma(\xi))}\leq C\Vert \widehat{f}(x,\eta) \Vert_{L^p(\mathbb{R}^2_x)}.
\end{equation}
Now, let us observe that
\begin{align*}
    \Vert \widehat{f}|_{\mathbb{T}^2}\Vert_{L^q({\mathbb{T}^2},d\sigma)} &=\Vert\widehat{f}(\xi,\eta) \Vert_{L^q((\mathbb{S}^1,d\sigma(\eta));L^q(\mathbb{S}^1,d\sigma(\xi)))} \\
    &\leq C\Vert \Vert \widehat{f}(x,\eta) \Vert_{L^p(\mathbb{R}^2_x)} \Vert_{L^q(\mathbb{S}^1,d\sigma(\eta))}=: C \Vert \widehat{f}(x,\eta)  \Vert_{L^q((\mathbb{S}^1,d\sigma(\eta)); L^p(\mathbb{R}^2_x) )}\\
    &:=I.
\end{align*}
Now, we will estimate the right hand side of the previous inequality. First, if we assume that $4/3\leq q<p'/3,$ then $p\leq q$ and the Minkowski integral inequality gives,
\begin{align*}
    I&=\left(\int\limits_{\mathbb{S}^1}\left(\int\limits_{\mathbb{R}^2}|\widehat{f}(x,\eta)|^pdx\right)^{\frac{q}{p}} d\sigma(\eta )\right)^{\frac{1}{q}} \leq \left(\int\limits_{\mathbb{R}^2}\left(\int\limits_{\mathbb{S}^1}|\widehat{f}(x,\eta)|^q d\sigma(\eta)\right)^\frac{p}{q} dx\right)^{\frac{1}{p}}\\
   & \lesssim \left(\int\limits_{\mathbb{R}^2}   \int\limits_{\mathbb{R}^2}|{f}(x,y)|^p dy dx\right)^{\frac{1}{q}}=\Vert f\Vert_{L^{p}(\mathbb{R}^4)},
\end{align*} where in the last inequality we have used the Fefferman restriction theorem. So we have proved that  \eqref{cardonaRest2018} holds true for $4/3\leq q<p'/3.$ Now, if $q<\frac{4}{3},$ then we  can use the finiteness of the measure $d\sigma(\xi,\eta)$ to deduce that
\begin{equation}
   \Vert \widehat{f}|_{\mathbb{T}^2}\Vert_{L^q({\mathbb{T}^2},d\sigma)}\lesssim \Vert \widehat{f}|_{\mathbb{T}^2}\Vert_{L^\frac{4}{3}({\mathbb{T}^2},d\sigma)} \leq C\Vert f\Vert_{L^{p}(\mathbb{R}^4)}
\end{equation} holds true for $1\leq p<\frac{4}{3}.$ Now, we will prove the converse announcement. So, let us assume that $p$ and $q$ are Lebesgue exponents satisfying \eqref{cardonaRest2018} with a constant $C>0$ independent of $f\in L^p(\mathbb{R}^4).$ If $g\in C^\infty_{c}(\mathbb{R}^2),$ let us define the function $f$ by $f(x,y)=g(x)g(y).$
The inequality,
\begin{equation}
   \Vert \widehat{f}|_{\mathbb{T}^2}\Vert_{L^q({\mathbb{T}^2},d\sigma)}:= \left(\int\limits_{\mathbb{T}^2}|\widehat{f}(\xi_1,\xi_2,\eta_1,\eta_2)|^qd\sigma(\xi_1,\xi_2,\eta_1,\eta_2)\right)^{\frac{1}{q}}\leq C\Vert f\Vert_{L^{p}(\mathbb{R}^4)},
\end{equation} implies that
\begin{equation}
   \Vert \widehat{g}|_{\mathbb{S}^1}\Vert_{L^q({\mathbb{S}^1},d\sigma)}:= \left(\int\limits_{\mathbb{S}^1}|\widehat{g}(\xi_1,\xi_2)|^qd\sigma(\xi_1,\xi_2)\right)^{\frac{1}{q}}\leq C\Vert g\Vert_{L^{p}(\mathbb{R}^2)}.
\end{equation}But, according with the Fefferman restriction theorem, the previous inequality only is possible for arbitrary $g\in C^\infty_{c}(\mathbb{R}^2),$ if $1\leq p<\frac{4}{3}$ and $q\leq p'/3.$

\section{Proof of Theorem \ref{ThrcardonaRestTn2018}}

Let us consider the $n$-dimensional torus
\begin{equation}
    \mathbb{T}^n=\{(x_{1,1},x_{1,2},x_{2,1},x_{2,2},\cdots,x_{n,1},x_{n,2}):x_{\ell,1}^2+x_{\ell,2}^2=1,\,1\leq \ell\leq n\}.
\end{equation}
We endow to $\mathbb{T}^n$ with the surface measure 
\begin{equation}
    d\sigma_n (\xi_{1,1},\xi_{1,2},\xi_{2,1},\xi_{2,2},\cdots,\xi_{n,1},\xi_{n,2})=\bigotimes_{j=1}^n d\sigma(\xi_{j,1},\xi_{j,2}),
\end{equation}
where $d\sigma$ is the `surface measure' on $\mathbb{S}^1.$ In order to prove Theorem \ref{ThrcardonaRestTn2018} we will use induction on $n.$ The case $n=2$ is precisely Theorem \ref{ThrcardonaRest2018}. So, let us assume that for some $n\in\mathbb{N},$ there exists $C_n$ depending only on the dimension $n,$ such that
\begin{equation}
   \Vert (\mathscr{F}_{\mathbb{R}^n}{u})|_{\mathbb{T}^{n}}\Vert_{L^q({\mathbb{T}^{n}},d\sigma_n)}\leq C_n\Vert u\Vert_{L^{p}(\mathbb{R}^{2n})},
\end{equation} 
for every function $u\in L^{p}(\mathbb{R}^{2n}).$ If $f\in  C^\infty_{c}(\mathbb{R}^{2n+2})\subset  L^p(\mathbb{R}^{2n+2}),$   $1\leq p<\frac{4}{3}$ and $q\leq p'/3,$  by using the approach of the previous section, we can write
\begin{align*}
    \widehat{f}(\xi_1,\xi_2,\eta)
    =\int\limits_{\mathbb{R}^2}e^{-i2\pi x\cdot \xi}(\mathscr{F}_{y\rightarrow \eta}{f}(x,\cdot))(\eta)dx,\,\,\eta\in \mathbb{R}^n.
\end{align*} By applying the Fefferman restriction theorem we deduce
\begin{equation}
    \Vert\widehat{f}(\cdot,\cdot,\eta) \Vert_{L^q(\mathbb{S}^1,d\sigma(\xi))}\leq \Vert \mathscr{F}_{y\rightarrow \eta}{f}(x,\cdot))(\eta) \Vert_{L^p(\mathbb{R}^2_x)}.
\end{equation} Now, by using that

\begin{align*}
 \Vert \widehat{f}|_{\mathbb{T}^{n+1}}\Vert_{L^q({\mathbb{T}^{n+1}},d\sigma_{n+1})} &=\Vert\widehat{f}(\xi_1,\xi_2,\eta) \Vert_{L^q((\mathbb{T}^n,d\sigma_n(\eta));L^q(\mathbb{S}^1,d\sigma(\xi)))} \\
    &\leq C\Vert \Vert \widehat{f}(x_{1,1},x_{1,2},\eta) \Vert_{L^p(\mathbb{R}^2_x)} \Vert_{L^q(\mathbb{T}^n,d\sigma_n(\eta))}\\
    &=: C \Vert \widehat{f}(x,\eta)  \Vert_{L^q((\mathbb{T}^n,d\sigma_n(\eta)); L^p(\mathbb{R}^2_x) )}\\
    &:=II,
    \end{align*} for $4/3\leq q<p'/3,$  $p\leq q,$ and  the Minkowski integral inequality, we have
\begin{align*}
    II&=\left(\int\limits_{\mathbb{T}^n}\left(\int\limits_{\mathbb{R}^2}|\widehat{f}(x,\eta)|^pdx\right)^{\frac{q}{p}} d\sigma_n(\eta )\right)^{\frac{1}{q}} \leq \left(\int\limits_{\mathbb{R}^2}\left(\int\limits_{\mathbb{T}^n}|\widehat{f}(x,\eta)|^q d\sigma_n(\eta)\right)^\frac{p}{q} dx\right)^{\frac{1}{p}}\\
   & \lesssim_n \left(\int\limits_{\mathbb{R}^{2}}   \int\limits_{\mathbb{R}^{2n}}|{f}(x,y)|^p dy dx\right)^{\frac{1}{p}}=\Vert f\Vert_{L^{p}(\mathbb{R}^{2n+2})},
\end{align*} where in the last inequality we have used the induction hypothesis. So, we have proved  Theorem \ref{ThrcardonaRestTn2018} for $4/3\leq q<p'/3.$ The case $q<\frac{4}{3}$ now follows from the finiteness of the measure $d\sigma_{n+1}$. That $1\leq p<\frac{4}{3}$ and $q\leq p'/3,$ are necessary conditions for \eqref{cardonaRestTn2018} can be proved if we replace $f$ in \eqref{cardonaRestTn2018} by a function of the form
\begin{equation}
    f(x_{1,1},x_{1,2},x_{2,1},x_{2,2},\cdots,x_{n,1},x_{n,2})=\prod_{j=1}^n g(x_{j,1},x_{j,2}),\,\,g \in C^\infty_c(\mathbb{R}^2).
\end{equation} Indeed, we automatically have \begin{equation}
   \Vert \widehat{g}|_{\mathbb{S}^1}\Vert_{L^q({\mathbb{S}^1},d\sigma)}:= \left(\int\limits_{\mathbb{S}^1}|\widehat{g}(\xi_1,\xi_2)|^qd\sigma(\xi_1,\xi_2)\right)^{\frac{1}{q}}\leq C\Vert g\Vert_{L^{p}(\mathbb{R}^2)}.
\end{equation} Consequently, the Fefferman restriction theorem, shows that the previous inequality only is possible for arbitrary $g\in C^\infty_{c}(\mathbb{R}^2),$ if $1\leq p<\frac{4}{3}$ and $q\leq p'/3.$\\
An usual argument of duality applied to Theorem \ref{ThrcardonaRestTn2018}, allows us to deduce the following result.
\begin{corollary}\label{ThrcardonaRestTn2018dual}
Let $F\in L^{q'}(\mathbb{T}^{n},d\sigma_n).$ Then there exists $C>0,$ independent of $F$ and satisfying
\begin{equation}\label{cardonaRestTn2018dual}
   \Vert (Fd\sigma_n)^{\vee} \Vert_{L^{p'}(\mathbb{R}^{2n})}=\left\Vert \int\limits_{ \mathbb{T}^n } e^{i2\pi x\cdot \xi}F(\xi)  d\sigma_n(\xi)   \right\Vert_{L^{p'}({\mathbb{R}^{2n}})}\leq C_n\Vert F\Vert_{L^{q'}{(\mathbb{T}^n,d\sigma_n)}},
\end{equation} if and only if $p'>4$ and $q'\geq (p'/3)'.$  Here, $d\sigma_n$ is the usual surface measure associated to $\mathbb{T}^{n}$ and $r':=r/r-1.$
\end{corollary}

\noindent {\bf Acknowledgement}. I would like to thanks Felipe Ponce from \textit{Universidad Nacional de Colombia} who introduced me to the restriction problem.

\bibliographystyle{amsplain}

\end{document}